\documentclass{amsart}
\usepackage{fancybox,color,graphicx,amsfonts,amsmath,amssymb,color}
\usepackage[all,color,line]{xy}
\newtheorem{theoremA}{Theorem}

\newtheorem{theorem}{Theorem}[section]

\theoremstyle{definition}

\theoremstyle{remark}
\newtheorem{remark}{Remark}[section]
\newtheorem{example}{Example}[section]

\numberwithin{equation}{section}

\newcommand{\mx}{\mathbf x}
\newcommand{\mX}{\mathbf X}

\newcommand{\mZ}{\mathbf Z}

\newcommand{\mA}{\mathbf A}

\newcommand{\mW}{\mathbf W}

\newcommand{\la}{\lambda}
\newcommand{\tr}{{\rm tr}}

\newcommand{\RR}{\mathbb{R}}
\newcommand{\CC}{\mathbb{C}}

\newcommand{\HH}{\mathbb{H}}

\newcommand{\E}{\mathbb{E}}

\newcommand{\calB}{\mathcal{B}}
\newcommand{\calM}{\mathcal{M}}

\author{W\l odzimierz Bryc}
\thanks{
\noindent Research partially supported by NSF
grant \#DMS-0504198 and Taft Research Seminar 2008/09.}
\address{
Department of Mathematical Sciences, University of Cincinnati, 2855
Campus Way, PO Box 210025, Cincinnati, OH 45221-0025, USA.}
\email{Wlodzimierz.Bryc@UC.edu}

\author{Virgil Pierce}
\thanks{Research partially supported by an NSF-VIGRE grant  \#DMS-0135308 and NSF grant \#DMS-0806219}
\address{Department of Mathematics, UT -- Pan American, 1201 W University Dr., Edinburg, TX  78539, USA.
}
\email{piercevu@utpa.edu}
\title{Duality of real and quaternionic random matrices}
\date{May 9, 2008; revised: January 14, 2009}
\begin{document}

\begin{abstract}
We show that quaternionic Gaussian  random
variables satisfy a generalization of the Wick formula for
computing the expected value of products in terms of a family of graphical
enumeration problems.  When applied to the quaternionic
Wigner and Wishart families of random matrices the result gives
the duality between moments of these families and the corresponding real Wigner
and Wishart families.
\end{abstract}
\maketitle

\section{Introduction}
The duality between symplectic and orthogonal groups has a long standing history,
and has been noted in physics literature in various settings, see e.g.
 \cite{Jungling-Oppermann-80,Mkrtchyan-81,Wegner-83,Witten-98}. Informally, the duality asserts that
 averages such as moments or partition functions
for the symplectic case of ``dimension'' $N$, can be derived from the respective
formulas  for the orthogonal case of dimension $N$ by inserting $-N$
 into these expressions and by simple scaling.
The detailed study of  the moments of one-matrix Wishart ensembles, with duality explicitly noted,
appears in \cite{Hanlon-Stanley-Stembridge-92}, see \cite[Corollary 4.2]{Hanlon-Stanley-Stembridge-92}.
The duality for one matrix Gaussian Symplectic Ensemble was noted
by Mulase and Waldron \cite{Mulase-Waldron-03} who introduced M\"obius graphs to write the
expansion for traces of powers of
GOE/GUE/GSE expansions in a unified way. The duality appears also in
 \cite[Theorem 6]{Ledoux-07} as a by-product of differential equations
 for the generating functions of moments.
Ref. \cite{Goulden-Jackson-96,Goulden-Jackson-96b,Goulden-Jackson-97}
analyze the related ``genus series" over locally orientable surfaces.

The purpose of this paper is to prove that the duality between moments of the Gaussian Symplectic Ensemble
and the Gaussian Orthogonal Ensemble, and between real Wishart and quaternionic Wishart ensembles extends
to several independent
matrices. Our technique consists of  elementary combinatorics; our proofs
 differ from \cite{Mulase-Waldron-03} in the one matrix case,
 and provide a more geometric interpretation for the duality; in the one-matrix Wishart case,
 our proof completes the combinatorial approach initiated in
  \cite[Section 6]{Hanlon-Stanley-Stembridge-92}.
The technique limits the scope of our results to moments, but the
relations between moments suggest similar relations between other analytic objects,
such as partition functions, see \cite{Mulase-Waldron-03}, \cite{Kodama-Pierce}.   The asymptotic expansion of the partition function and analytic description of
 the coefficients of this expansion
  for $\beta=2$ case appear in \cite{Ercolani-McLaughlin-03,Guionnet-Maurel-Segala-0503064,Maurel-Segala-0608192}.

The paper is organized as follows. In Section \ref{Sect1} we review basic properties of quaternionic
Gaussian random variables. In Section \ref{Sect2} we introduce M\"obius graphs;  Theorems \ref{T quaternion moments} and \ref{thm2.1}
 give formulae
for the expected values of products of quaternionic Gaussian random
variables in terms of the Euler characteristics of sub-families of
M\"obius graphs or of bipartite M\"obius graphs.
In Section \ref{duality} we apply the formulae to the quaternionic Wigner and  Wishart families.

In this paper, we do not address the question of whether the duality
 can be extended to more general functions, or to more general
 $\beta$-Hermite and $\beta$-Laguerre ensembles introduced in \cite{Dumitru-Edelman-06}.

\section{ Moments of quaternion-valued Gaussian random variables}\label{Sect1}

\subsection{Quaternion Gaussian law}
Recall that a quaternion $q\in \HH$ can be represented as $q=x_0+i x_1+j x_2+ k x_3$
with
$i^2=j^2=k^2=ijk=-1$ and with real coefficients $x_0,\dots,x_3$.
The conjugate quaternion is $\overline{q}=x_0-i x_1-j x_2- k x_3$,
so $|q|^2:=q\overline{q}\geq 0$. Quaternions with $x_1=x_2=x_3=0$ are usually identified with real numbers;
the real part of a quaternion is  $\Re(q)=(q+\bar{q})/2$.

It is well known that quaternions can be identified with the set of certain $2\times2$ complex matrices:
\begin{equation}
  \label{H2C}
 \HH\ni x_0+ix_1+jx_2+kx_3\sim \left[\begin{matrix}
  x_0+ix_1 & x_2+i x_3 \\
  -x_2+ix_3&x_0-ix_1
\end{matrix}\right]\in \calM_{2\times 2}(\CC),
\end{equation}
where on the right hand side $i$ is the usual imaginary unit of $\CC$.
Note that since $\Re(q)$ is twice the trace of the matrix
representation in \eqref{H2C}, this implies the cyclic property
\begin{equation}
  \label{tmp**}\Re(q_1q_2\dots q_n)=\Re(q_2q_3\dots q_nq_1).
\end{equation}

The (standard) quaternion Gaussian random variable is an $\HH$-valued random variable
which can be represented as
\begin{equation}
  \label{HH}
  Z=\xi_0+i \xi_1+ j\xi_2+k\xi_3
\end{equation} with independent real
normal $N(0,1)$ random variables $\xi_0,\xi_1,\xi_2,\xi_3$.
Due to symmetry of the centered normal laws on $\RR$,
the law of $(Z,\overline{Z})$ is the same as the law of $(\overline{Z},Z)$.
A  calculation shows  that if $Z$ is quaternion Gaussian then
for fixed $q_1, q_2 \in\HH$,
$$
\E(Z q_1 Zq_2)=\E(Z^2)\bar{q}_1q_2,\;
\E(Z q_1 \overline{Z}q_2)=\E(Z\bar{Z})\Re(q_1)q_2\,.
$$
For future reference, we insert explicitly the moments:
\begin{eqnarray}
 \label{q1}
\E(Z q_1 Zq_2)&=&-2\bar{q}_1q_2,\\
\label{q2}
\E(Z q_1 \overline{Z}q_2)&=& 2 (q_1+\bar{q}_1)q_2.
\end{eqnarray}
By linearity, these formulas imply
\begin{eqnarray}
\label{q3}
\E( \Re(Z q_1) \Re( \overline{Z} q_2) ) &=& \Re( q_1 q_2 ) ,\\
\label{q4}
\E( \Re(Z q_1) \Re(Z q_2 ) ) &=& \Re( \bar{q}_1 q_2 ).
\end{eqnarray}

\subsection{Moments}
The following is  known as Wick's theorem  \cite{Wick-50}.
\begin{theoremA}[Isserlis \cite{Isserlis-1918}]
\label{WickTHMR} If $(X_1,\dots,X_{2n})$ is
a $\RR^{2n}$-valued Gaussian random vector with mean zero, then
\begin{equation}\label{R-Wick}
  E(X_1X_2\dots X_{2n})=\sum_{V}\prod_{\{j,k\}\in V} E(X_jX_k),
\end{equation}
where the sum is taken over all pair partitions  $V$ of
$\{1,2,\dots,2n\}$, i.e.,
partitions into two-element sets, so each $V$ has the form
$$V=\left\{\{j_1,k_1\},\{j_2,k_2\},\dots,\{j_n, k_n\}\right\}.$$
\end{theoremA}

%
%

 Theorem \ref{WickTHMR} is a consequence of the moments-cumulants relation \cite{Leonov-Shirjaev-59};
  the connection is best visible in the partition formulation of \cite{Speed-83}. For another proof,
  see \cite[page 12]{Janson-97}.

  Our first goal is to extend this formula to certain quaternion Gaussian random variables.
The  general multivariate quaternion Gaussian law is discussed in
\cite{Vakhania-99}. Here we will only consider a special setting of
sequences that are drawn with repetition from a sequence of
independent standard Gaussian quaternion random variables. In
section \ref{duality}
 we apply this result to a multi-matrix version of the
duality between GOE and GSE ensembles of random matrices.

In view of the Wick formula \eqref{R-Wick} for real-valued jointly Gaussian random
variables, formulas \eqref{q1} and \eqref{q2} allow us to compute
moments of certain products of quaternion Gaussian random variables.
Suppose  the $n$-tuple $(X_1,X_2,\dots,X_{n})$
consists of random variables taken, possibly with repetition, from
the set
 $$\{Z_1,\bar{Z}_1,Z_2,\bar{Z_2},\dots\},$$
 where
$Z_1,Z_2,\dots$ are independent quaternion Gaussian random variables.
 Consider an auxiliary family of independent pairs $\{(Y_{j}^{(r)},Y_{k}^{(r)}): r=1,2,\dots\}$ which have the same laws as
 $(X_j,X_k)$, $1\leq j,k\leq n$ and are independent for different $r$.
  Then the Wick formula for real-valued Gaussian variables implies $\E(X_1X_2\dots X_{n})=0$ for odd $n$, and
   \begin{equation}\label{Wick0}
\E(X_1X_2\dots X_{n})=\sum_{f}\E(Y_1^{(f(1))}Y_2^{(f(2))}\dots Y_{n}^{(f(n))}),
 \end{equation}
 where the sum is over the pair partitions $V$  that appear under the sum in Theorem \ref{WickTHMR}, each represented by the level sets of a two-to-one valued function
 $f:\{1,\dots,n\}\to\{1,\dots,m\}$ for $n=2m$.  (Thus the sum is over classes of equivalence of $f$, each of $m!$ representatives contributing the same value.)

 For example, if $Z$ is quaternion Gaussian then
applying \eqref{Wick0} with $f_1$ that is constant, say $1$,  on $\{1,2\}$,
  $f_2$ that is constant  on $\{1,3\}$, and $f_3$ that is constant   on $\{1,4\}$
we get
$$
\E(Z^4)=\E(Z^2) \left(\E(Z^2) + \E(\bar{Z}Z)+ \E(\bar{Z}^2)\right)=0.
$$

Formulas \eqref{q1} and \eqref{q2} then show that the Wick reduction step takes the
following form.
\begin{equation}\label{Wick1}
  \E(X_1X_2\dots X_n)=\sum_{j=2}^n \E(X_1X_j)\E(U_j X_{j+1}\dots X_n),
\end{equation}
where
$$
U_j=\begin{cases}
  \Re(X_2\dots X_{j-1}) & \mbox{ if $X_j=\bar{X}_1$ }\\
  \bar{X}_{j-1}\dots \bar{X}_2 &\mbox{ if $X_j={X}_1$}\\
  0 &\mbox{ otherwise} \,.
\end{cases}
$$
This implies that one can do inductively the calculations, but due to noncommutativity
arriving at explicit answers may still require significant work.

Formula \eqref{Wick1} implies that $\E(X_1X_2\dots X_n)$ is real, so on the left hand side of \eqref{Wick1} we can write
$\E(\Re(X_1X_2\dots X_n))$; this form of the formula will be associated with one-vertex M\"obius graphs.

Furthermore, we have a Wick reduction which will correspond to the multiple vertex M\"obius graphs:
\begin{multline}\label{Wick2}
\E( \Re( X_1 ) \Re( X_2 X_3 \dots X_n ) ) \\=
\sum_{j=2}^n \E(\Re(X_1)\Re(X_j))\E( \Re(X_2 \dots X_{j-1} X_{j+1} \dots X_n  )).
\end{multline}
(This is just a consequence of Theorem \ref{WickTHMR}).

In the next section we will show that formulae (\ref{Wick1}) and
(\ref{Wick2}) give a method of computing the expected values of
quaternionic Gaussian random variables by the enumeration of M\"obius
graphs partitioned by their Euler characteristic.  This is analogous
to similar results for complex Gaussian random variables and ribbon
graphs, and for real Gaussian random variables and M\"obius graphs.
In Section \ref{duality} we will show  that this
result implies the duality of the GOE and GSE
ensembles of Wigner random matrices and the duality of real and
quaternionic Wishart random matrices.

\section{M\"obius graphs and quaternionic Gaussian moments}\label{Sect2}

In this section we introduce M\"obius graphs and then give formulae
for the expected values of products of quaternionic Gaussian random
variables in terms of the Euler characteristics of sub-families of
M\"obius graphs.  This is an analogue of the method of t'Hooft
\cite{Bessis-Itzykson-Zuber-80,Goulden-Harer-Jackson-01, Goulden-Jackson-97,harer-zagier-86, Jackson-94, tHooft}.  M\"obius graphs have also been used to give combinatoric
interpretations of the expected values of traces of Gaussian
orthogonal ensemble of random matrices and of Gaussian symplectic ensembles, see the articles
\cite{Goulden-Jackson-97} and \cite{Mulase-Waldron-03}.  The  connection between M\"obius graphs and
quaternionic Gaussian random variables is at the center of the work of Mulase
and Waldron \cite{Mulase-Waldron-03}.

\subsection{M\"obius graphs}
M\"obius graphs are ribbon graphs where the edges (ribbons) are
allowed to twist, that is they either preserve or reverse the local
orientations of the vertices.  As the convention is that the ribbons
in \textit{ribbon graphs} are not twisted we follow \cite{Mulase-Waldron-03} and call the unoriented
variety M\"obius graphs. The vertices of a M\"obius graph are
represented as disks together with a local orientation; the edges are represented as
ribbons, which preserve or reverse the local orientations of the
vertices connected by that edge. Next we identify the collection of
disjoint cycles of the sides of the ribbons found by following the
sides of the ribbon and obeying the local orientations at each
vertex.  We then attach disks to each of these cycles by gluing the
boundaries of the disks to the sides of the ribbons in the cycle.
These disks are called the faces of the M\"obius graph, and the
resulting surface we find is the surface of maximal
Euler
characteristic on which the M\"obius graph may be drawn so that edges do not cross.

Denote by $v(\Gamma)$, $e(\Gamma)$, and $f(\Gamma)$ the number of vertices, edges, and faces of $\Gamma$.
We say that the Euler characteristic of $\Gamma$ is
$$
\chi(\Gamma)=v(\Gamma)-e(\Gamma)+f(\Gamma),
$$
for connected $\Gamma$, this is also the maximal Euler characteristic of a connected surface into which
$\Gamma$ is embedded.
For example, in Fig. \ref{F1}, the Euler characteristics are
$\chi_1=1-1+2=2$ and $\chi_2=1-1+1=1$.  The two graphs may be embedded
into the sphere or projective sphere respectively.

If  $\Gamma$ decomposes into connected components $\Gamma_1$, $\Gamma_2$, then
$\chi(\Gamma)=\chi(\Gamma_1)+\chi(\Gamma_2)$.

Throughout the paper our M\"obius graphs will have the following
labels attached to them:  the vertices are labeled to make them
distinct, in addition the edges emanating from each vertex are also
labeled so that rotating any vertex produces a distinct graph. These
labels may be removed if one wishes by rescaling all of our
quantities by the number of automorphisms the unlabeled graph would
have.

\subsection{Quaternion version of Wick's theorem}

\tolerance=2000
Suppose  the $2n$-tuple $$(X_1,X_2,\dots,X_{2n})$$
consists of random variables taken, possibly with repetition, from
the set $\{Z_1,\bar{Z_1},Z_2,\bar{Z_2},\dots\}$ where
$Z_1,Z_2,\dots$ are independent quaternionic Gaussian.
Fix a sequence $j_1,j_2,\dots,j_m$ of natural numbers such that $j_1+\dots+j_m=2n$.

Consider the family $\calM=\calM_{j_1,\dots,j_m}(X_1,X_2,\dots,X_{2n})$, possibly empty, of M\"obius graphs with $m$ vertices of degrees $j_1,j_2,\dots,j_m$ with edges labeled by $X_1,X_2,\dots,X_{2n}$,
whose regular  edges  correspond to pairs $X_i=\bar{X}_j$ and flipped edges  correspond to pairs $X_i=X_j$. No edges of  $\Gamma\in \calM$ can join random variables $X_i,X_j$  that are independent.
\tolerance=1000

\begin{theorem}
  \label{T quaternion moments}
Let $\left\{ X_1, X_2, \dots, X_{2n} \right\}$ be chosen, possibly
with repetition, from the set $\{ Z_1, \bar{Z_1}, Z_2, \bar{Z_2},
\dots \}$ where $Z_j$ are independent quaternionic Gaussian random
variables, then
\begin{multline}
\E\big(\Re(X_1X_2 \dots X_{j_1})\Re(X_{j_1+1}^{}\dots X_{j_1+j_2}^{})\times\dots \\\times\Re(X_{j_1+j_2+\ldots+j_{m-1}+1}^{}\dots
  X_{2n}^{})\big)=4^{n-m}\sum_{\Gamma\in \calM} (-2)^{\chi(\Gamma)}.
\end{multline}
(The right hand side is interpreted as $0$ when $\calM=\emptyset$.)
\end{theorem}

\begin{remark}
  We would like to emphasize that in computing the Euler characteristic one must first break the graph into
  connected components. For example, if $j_1=\dots=j_m=1$ so that $m=2n$ is even,
and $X_1, X_2, \dots, X_{2n}$ are $n$ independent pairs, as the real
parts are commutative we may assume that $X_{2k} = X_{2k-1}$, and
the moment is
$$ 1=\E( \Re (X_1) \Re (X_2) \dots \Re (X_{2n}) )  = 4^{-n}  (-2)^{\chi(\Gamma)}.  $$
We see that graphically $\Gamma$
 is a collection of $2n$ degree one vertices connected together forming $n$
 dipoles (an edge with a vertex on either end).
  Hence there are $n$ connected components each of Euler characteristic $2$,
  therefore the total Euler characteristic is $\chi = 2 n $  giving
$$ 4^{-n}  (-2)^\chi= 4^{-n} 4^{n} = 1. $$
\end{remark}

\begin{proof}[Proof of Theorem \ref{T quaternion moments}]
In view of \eqref{Wick0} and \eqref{Wick1}, it suffices to show that  if $X_1,\dots,X_{2n}$ consists of $n$ independent pairs, and each pair is either of the form $(X,X)$ or $(X,\bar{X})$, then
\begin{multline}
\label{star2}
  \E\big(\Re(X_1X_2 \dots X_{j_1})\Re(X_{j_1+1}^{}\dots X_{j_1+j_2}^{})\times\dots \\\times\Re(X_{j_1+j_2+ \dots+ 
  j_{m-1}+1}^{}\dots
  X_{2n}^{})\big)=  4^{n-m}(-2)^{\chi(\Gamma)}\;,
\end{multline}
where $\Gamma$ is the M\"obius graph that describes the pairings of the sequence.

First we check the two M\"obius graphs for
$n=1$, $m=1$:
\[   \E( \Re( X \bar{X}) ) = (-2)^2  \,, \quad \mbox{and}\quad \E( \Re(X X))
= (-2)^1   \,.\]
One checks that these correspond to the M\"obius graphs in Figure \ref{F1},
which gives a sphere ($\chi=2$)  and projective sphere ($\chi=1$) respectively.
\begin{figure}[htb]
 \includegraphics[width=4in]{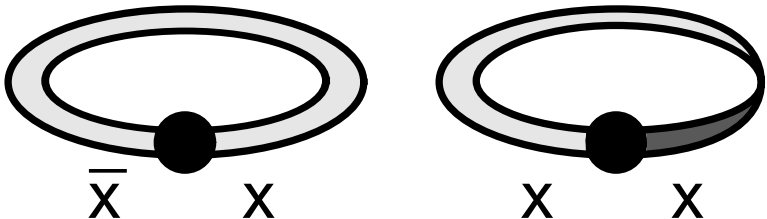}
 \caption{The two possible M\"obius graphs with a single degree 2 vertex. The left hand one is a ribbon which is untwisted and the graph embeds into a copy of the Riemann sphere, while the
 right hand one is a ribbon which is twisted and the graph embeds into a copy of the projective sphere.  
 \label{F1}}
\end{figure}

We now proceed with the induction step.
 One notes that by independence of the pairs at different edges, the
left hand side of \eqref{star2} factors into the product corresponding to
connected components of $\Gamma$. It is therefore enough to consider
connected $\Gamma$. \label{tmp1}

If $\Gamma$ has two vertices that are joined by an edge, we can use cyclicity
of $\Re$  to move the variables that label the edge to the first positions in
their cycles, say $X_1$ and $X_{j_1+1}$ and use \eqref{q3} or \eqref{q4} to
eliminate this pair from the product. The use of relation (\ref{q3})
is just that of gluing the two vertices together removing the edge $x$
which is labeled by the two appearances of $Z$.  Relation (\ref{q4}) glues
together the two vertices, removing
the edge $x$, and the reversal of orientation across the edge is
given by the conjugate (see Figure \ref{F4}).
These geometric operations reduce $n$ and $m$ by one without changing the
Euler characteristic:
the number of edges and the number of vertices are reduced by 1; the faces are
preserved -- in the case of edge flip
in Fig. \ref{F4}, the edges of the face from which we remove the edge, after
reduction follow the same order.

 Therefore we will only need to prove the
result for the single vertex case of the induction step.

\begin{figure}[htb]
\includegraphics[height=4in]{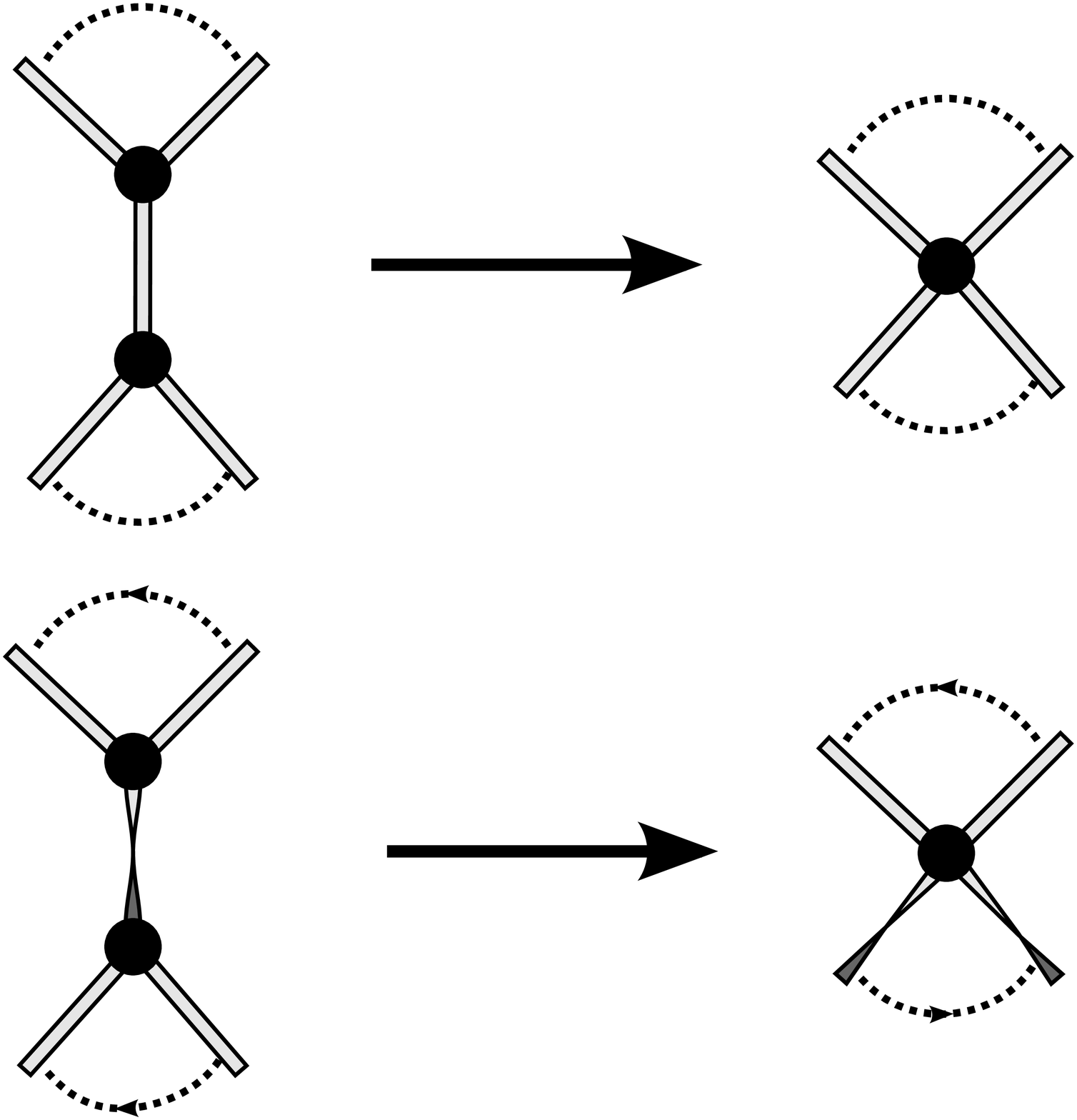}
  \caption{A M\"obius graph with two vertices connected by a ribbon may be reduced to a M\"obius graph with one less vertex and one less edge.
In these graphs the ``$\dots$" are to mean that there are arbitrary
numbers of other edges at the vertex, and the edges drawn are
connected to other vertices.
 The top graph is an example of this reduction when the connecting ribbon is untwisted, in this case the two vertices are glued together with no other changes in the ribbons.  The bottom graph is an example of this reduction when the connecting ribbon is twisted, in this case the two vertices are glued together and the order and orientation (twisted or untwisted) of the ribbons on one side are reversed.
\label{F4}}
\end{figure}

We wish to show that
\begin{equation} \label{star3}
 \E( X_1 X_2 X_3 X_4 \dots X_{2n}) =
(-2)^{\chi({\Gamma})}  4^{n-1} \,,
\end{equation}
where $\Gamma$ is a one vertex M\"obius graph  with arrows (half edges)
 labeled by
$X_k$.
We will do this by induction, there are two cases:
\begin{enumerate}

\item[\bf Case 1:] $X_1 = \bar{X}_j$ for $1< j \leq 2n $,
\begin{align} \nonumber
\E( X_1 X_2 \dots X_{j-1} \bar{X}_1 X_{j+1} \dots X_{2n}) &=
\E( X_1 \bar{X}_1 ) \E( \Re( X_2 \dots X_{j-1} ) \Re( X_{j+1} \dots X_{2n}) ) \\
\label{oriented_reduction}
&= 4 \E( \Re( X_2 \dots X_{j-1}) \Re( X_{j+1} \dots X_{2n}) ) \,.
\end{align}
This corresponds to the reduction of the M\"obius graph pictured
in Figure \ref{orient-reduct}, which splits the single vertex
into two vertices.  The Euler characteristic becomes $\chi_2 = 2 - (n-1) + f_1
= \chi_1 +2 $.
By the induction assumption we find
\begin{equation*}
\E( X_1 \dots X_{2n}) = 4 \left[ 4^{(n-1) - 2} (-2)^{\chi_2} \right]
= 4^{n-2} (-2)^{\chi_1 +2} = 4^{n-1} (-2)^{\chi_1} \,.
\end{equation*}

\begin{figure}
\begin{center}
\includegraphics[width=8cm]{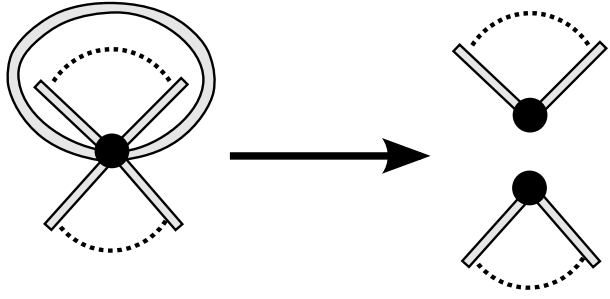}
\end{center}
\caption{ Here we have an untwisted ribbon of the M\"obius graph
returning to the same vertex, this edge is removed in our reduction
procedure giving us a M\"obius graph with one more vertex and one
less edge. \label{orient-reduct}}
\end{figure}

\item[\bf Case 2:] $X_1 = X_j$ for $1 < j \leq 2n$,
\begin{align}\nonumber
\E( X_1 X_2 \dots X_{j-1} X_1 X_{j+1} \dots X_{2n}) &=
\E( X_1 X_1 ) \E(  \bar{X}_{j-1} \dots \bar{X}_2  X_{j+1} \dots X_{2n} )
\\ \label{unoriented-reduction}
&= (-2) \E( \bar{X}_{j-1} \dots \bar{X}_2 X_{j+1} \dots X_{2n} ) \,.
\end{align}
This corresponds to the reduction of the M\"obius graph pictured
in Figure \ref{unorient-reduct}, which keeps the single vertex
and flips the order and orientation of the edges between $X_1$ and $X_j$.
The Euler characteristic becomes $\chi_2 = 1 - (n-1) + f_1 = \chi_1 + 1$.
By the induction assumption we find
\begin{equation*}
\E( X_1 \dots X_{2n}) = (-2) \left[ 4^{(n-1) -1} (-2)^{\chi_2} \right]
= (-2)^{-1} 4^{n-1} (-2)^{\chi_1 + 1} = 4^{n-1} (-2)^{\chi_1}\,.
\end{equation*}

\begin{figure}
\begin{center}
\includegraphics[width=8cm]{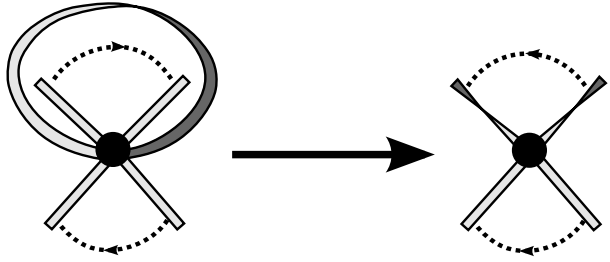}
\end{center}
\caption{ Here we have a twisted ribbon of the M\"obius graph
returning to the same vertex, this edge is removed in our reduction
procedure giving us a M\"obius graph with one less edge and giving a
reverse in both the order and orientation (twisted or untwisted) of
the ribbons on one side of the removed ribbon.
\label{unorient-reduct}}
\end{figure}

\end{enumerate}

Note that taking away an oriented ribbon creates a new vertex. The
remaining graph might still be connected, or it may split into two
components. If taking away a loop makes the graph disconnected, then
the counts of changes to edges and vertices are still the same. But
the faces need to be counted as follows: the inner face of the
removed edge becomes the outside face of one component, and the
outer face at the removed edge becomes the outer face of the other
component. Thus the counting of faces is not affected by
whether the graph is connected.

With these two cases checked, by the induction hypothesis, the
proof is completed.

\end{proof}

\subsection{Bipartite M\"obius graphs and quaternionic Gaussian moments} \label{Sect bipartite}
To deal with quaternionic Wishart random matrices, we need to consider a special subclass of quaternionic
Gaussian variables from Theorem \ref{T quaternion moments}.
Suppose  the $2n$-tuple $(X_{\pm 1},X_{\pm 2},\dots,X_{\pm n})$
consists of $n$ pairs of random variables taken with repetition, from
the set $\{Z_1,Z_2,\dots,Z_n\}$ of independent quaternionic Gaussian random variables. Note that in contrast to the setup for Theorem \ref{T quaternion moments}, here all $Z$'s are without conjugation.
Fix a sequence $j_1,j_2,\dots,j_m$ of natural numbers such that $j_1+\dots+j_m=n$.
Theorem \ref{T quaternion moments} then says that
\begin{align*}
 \E\big( \Re( \bar{X}_{-1} X_1 \cdots \bar{X}_{-j_1} X_{j_1})
\times \Re( \bar{X}_{-j_1-1} X_{j_1 +1} \cdots \bar{X}_{-j_1-j_2} X_{j_1+j_2} )
\times \cdots
\\ \dots \times
\Re( \bar{X}_{-j_1-j_2-\dots -j_{m-1}-1} \cdots \bar{X}_{-n} X_{n} )\big)
\\ = 4^{n-m} (-2)^{\chi(\Gamma)},
\end{align*}
where $\Gamma$ is the  M\"obius graph with $m$ vertices and with edges labeled by $X_{\pm 1},\dots,X_{\pm n}$ that describes the
pairings between the variables under the expectation. Our goal is to show that the same formula holds true for another graph, a bipartite M\"obius graph whose edges are labeled by $n$ pairs $(X_{-j},X_j)$, $1\leq j \leq n$.

The bipartite M\"obius  graph has two types of vertices: black vertices and white (or later, colored) vertices with ribbons that can only connect a black vertex to a white vertex. (As previously, the ribbons may carry a ``flip" of orientation which we represent graphically as a twist.) To define this graph, we need to introduce three pair partitions on the set $\{\pm 1,\dots,\pm n\}$.
The first partition, $\delta$, pairs $j$ with $-j$. The second partition, $\sigma$, describes the placement of $\Re$: its pairs are
\begin{align*}
\{1,-2\},\{2,-3\},\dots,\{j_1,-1\},\\
\{j_1+1, -(j_1+2)\},\dots, \{j_1+j_2,-(j_1+1)\},\\
\vdots\\
\{1+\sum_{k=1}^{m-1}j_k,-(2+\sum_{k=1}^{m-1}j_k)\},\dots, \{n,-(1+\sum_{k=1}^{m-1}j_k)\}.
\end{align*}
The third partition, $\gamma$, describes the choices of pairs from $Z_1,\dots,Z_n$. Thus
$
\{j,k\}\in\gamma
$
if $X_j=X_k$ when $jk>0$ or $X_j=\bar{X}_k$ if $jk<0$.

We will also represent these pair partitions as graphs with vertices arranged in
two rows, and
with the edges drawn between the vertices in each pair of  a partition.
Thus
\begin{equation}
  \label{eq:delta}
  \delta=\begin{matrix}
  \xymatrix  @-1pc{
 {^1_ \bullet} \ar@{-}[d]& {^2_\bullet}
  \ar@{-}[d]& \dots & {^n_\bullet} \ar@{-}[d]\\
 {^{\hspace{1.5mm}\bullet}_{-1}} &{^{\hspace{1.5mm}\bullet}_{-2}} &\dots &
{^{\hspace{1.5mm}\bullet}_{-n}} \\
}
\end{matrix}
\end{equation}
and
{\small
$$
\sigma=\begin{matrix}
 \xymatrix  @-1pc{
 {^1_ \bullet} \ar@{-}[dr]& {^2_\bullet} \ar@{-}[dr]& {^3_\bullet}
 & \dots & {^{j_1}_\bullet}\ar@{-}[dllll] & {^{j_1+1}_{\hspace{2.5mm}\bullet}}\ar@{-}[dr] & {^{j_1+2}_{\hspace{3mm}\bullet}}&\dots & {^{j_1+j_2}_{\hspace{3.7mm}\bullet}}\ar@{-}[dlll]&
  {^{j_1+j_2+1}_{\hspace{5mm}\bullet}}&\dots \\
 {^{\hspace{1.5mm}\bullet}_{-1}} &{^{\hspace{1.5mm}\bullet}_{-2}} &{^{\hspace{1.5mm}\bullet}_{-3}}&\dots &
{^{}_\bullet}& {^{}_\bullet} & {^{}_\bullet}&\dots&{^{}_\bullet}&{^{}_\bullet}&\dots
}
\end{matrix}
$$
}

Consider the 2-regular graphs $\delta\cup \gamma$ and $\delta\cup \sigma$. We orient the cycles of these graphs by ordering $(-j,j)$ on the left-most vertical edge of the cycle.
For example,

{\small $$
\sigma\cup\delta=\begin{matrix}
 \xymatrix  @-1pc{
 {^1_ \bullet} \ar@{-}[dr]& {^2_\bullet} \ar@{-}[dr]& {^3_\bullet}
 & \dots & {^{j_1}_\bullet}\ar@{-}[dllll] & {^{j_1+1}_{\hspace{2.5mm}\bullet}}\ar@{-}[dr] & {^{j_1+2}_{\hspace{3mm}\bullet}}&\dots & {^{j_1+j_2}_{\hspace{3.7mm}\bullet}}\ar@{-}[dlll]& {^{j_1+j_2+1}_{\hspace{5mm}\bullet}}&\dots \\
 {^{\hspace{1.5mm}\bullet}_{-1}}\ar@{->}[u] &{^{\hspace{1.5mm}\bullet}_{-2}}\ar@{-}[u] &{^{\hspace{1.5mm}\bullet}_{-3}}\ar@{-}[u]&\dots &
{^{}_\bullet}\ar@{-}[u]& {^{}_\bullet} \ar@{->}[u]& {^{}_\bullet}\ar@{-}[u]&\dots&{^{}_\bullet}\ar@{-}[u]& {^{}_\bullet} \ar@{->}[u]&\dots
}
\end{matrix}
$$
}
We now define the bipartite M\"obius graph by assigning black vertices to the $m$ cycles of  $\delta\cup \sigma$, and white vertices to the cycles of  $\delta\cup \gamma$.

Each black vertex is oriented counter-clockwise. For each black vertex we follow the cycle of $\delta\cup \sigma$, drawing a labeled line for each element of the partition. The lines corresponding to $-j,j$ are adjacent, and will eventually become two edges of a ribbon.

Each white vertex is oriented, say, clockwise; we identify the graphs that differ only by
a choice of orientation at some of the white vertices.
For each white vertex, we follow the corresponding cycle of $\delta\cup \gamma$, drawing a labeled line for each element of the partition. The lines corresponding to $-j,j$ are adjacent, but may appear in two different orders depending on the orientation of the corresponding edge of $\delta$ on the cycle.

The final step is to connect pairs $(j,-j)$ on the black vertices with the same pairs on the white vertices. This creates the ribbons, which carry a flip if the orientation of the two lines
on the black vertex does not match the orientation of the same edges on the white vertex.

\begin{figure}
\begin{center}\includegraphics[width=6cm]{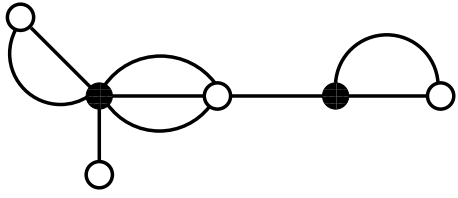}
\end{center}
\caption{\label{black_white} Representation of a bipartite M\"obius
graph, the edges drawn are ribbons that are either twisted or
untwisted.  }
\end{figure}
\begin{figure}
\begin{center}\includegraphics[width=5cm]{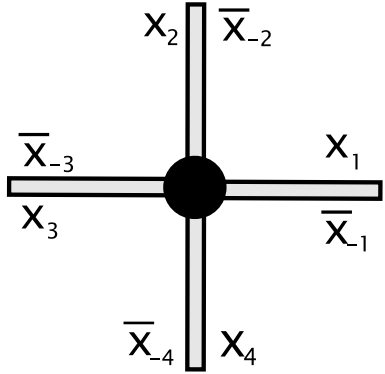}\end{center}
\caption{\label{black-white-2} Example of the labeling we use for
the edges emanating from a black vertex in a bipartite M\"obius
graph.  }
\end{figure}

The individual edges pictured in Figure \ref{black_white} are ribbons and are
labeled as in Figure \ref{black-white-2}.
  We allow twists of ribbons to
propagate through a white vertex calling the two bipartite M\"obius graphs
equivalent.

Suppose  the $2n$-tuple $(X_{\pm 1},X_{\pm 2},\dots,X_{\pm n})$
consists of random variables taken with repetition, from
the set $\{Z_1,Z_2,\dots,Z_n\}$ of independent quaternionic Gaussian random variables.
Let $\calM=\calM(X_{\pm 1},X_{\pm 2},\dots,X_{\pm n})$ denote the set of all bipartite M\"obius graphs $\Gamma$
that correspond to various ways of pairing all repeated $Z$'s in the sequence
$(X_{\pm 1},X_{\pm 2},\dots,X_{\pm n})$; the pairs are given by adjacent half edges at each white
vertex. (See the preceding construction.) $\calM=\emptyset$ if there is a $Z_j$ that is repeated
an odd number of times.
\begin{theorem} \label{thm2.1}
\begin{align*}
 \E\big( \Re( \bar{X}_{-1} X_1 \cdots \bar{X}_{-j_1} X_{j_1})
\times \Re( \bar{X}_{-j_1-1} X_{j_1 +1} \cdots \bar{X}_{-j_1-j_2} X_{j_1+j_2} )
\times \cdots
\\ \dots \times
\Re( \bar{X}_{-j_1-j_2-\dots -j_{m-1}-1} \cdots \bar{X}_{-n} X_{n} )\big)
\\ = 4^{n-m}\sum_{\Gamma\in\calM} (-2)^{\chi(\Gamma)}\;,
\end{align*}
(The right hand side is interpreted as $0$ when $\calM=\emptyset$.)
\end{theorem}

\begin{proof}
The proof is fundamentally the same as that of the Wigner version of this
theorem.
In view of \eqref{Wick0} and \eqref{Wick1}, it suffices
 to consider  $\{ X_{\pm 1}, \dots X_{\pm n} \}$  that  form $n$ independent pairs, and show that
\begin{multline}\label{WWW314}
 \E\big( \Re( \bar{X}_{-1} X_1 \cdots \bar{X}_{-j_1} X_{j_1})
\times \Re( \bar{X}_{-j_1-1} X_{j_1 +1} \cdots \bar{X}_{-j_1-j_2} X_{j_1+j_2} )
\times \cdots
\\  \dots \times
\Re( \bar{X}_{-j_1-j_2-\dots -j_{m-1}-1} \cdots \bar{X}_{-n} X_{n} )\big)
 = 4^{n-m} (-2)^{\chi(\Gamma)}\;,
\end{multline}
where $\Gamma$ is the bipartite M\"obius graph that describes the
pairings.

We will prove \eqref{WWW314} by induction; to that end we first check
that with $n=1, m=1$
we have $\E(\bar{X}X)=(-2)^2$ in agreement with \eqref{WWW314}.

If $\Gamma$ has two black vertices connected together by edges adjacent at a white vertex,
we can use the cyclicity of $\Re$ to move the variables that label the
respective edges and share the same face to
the first position in their cycles, so that we may call them $\bar{X}_{-1}$
and either $X_j$ or $\bar{X}_{-j}$.
We now use relations (\ref{q3}) and (\ref{q4}) to eliminate the pair from the
product:
\begin{equation*}
\E\left( \Re( \bar{X}_{-1} \cdots X_{j_1}) \Re( X_j \dots X_{j_1+j_2} \bar{X}_{-j} )\right)
=  \E( X_1 \cdots X_{j_1+j_2} \bar{X}_{-j})\,,
\end{equation*}
or
\begin{equation*}
\E\left( \Re( \bar{X}_{-1} \cdots X_{j_1}) \Re( \bar{X}_j \dots X_{j_1+j_2} )\right)
= \E( \bar{X}_{j_1} X_{-j_1} \cdots \bar{X}_1 X_j   \cdots X_{j_1+j_2})\,.
\end{equation*}
The use of relation (\ref{q3}) corresponds to that of gluing together the two
ribbons along the halves adjacent at the white vertex, and gluing together the
corresponding black vertices (see Figure \ref{two-gluing-orient}).
The use of relation (\ref{q4}) corresponds to the same gluing, but in this
case one of the ribbons has an orientation reversal in it, resulting in an
orientation reversal for the remaining sides (see Figure
\ref{two-gluing-unorient}).
These geometric operations reduce $n$ and $m$ by one without changing the
Euler characteristic: both the number of edges and the number of vertices are
reduced by one while the number of faces is preserved.

\begin{figure}
\includegraphics[width=12cm]{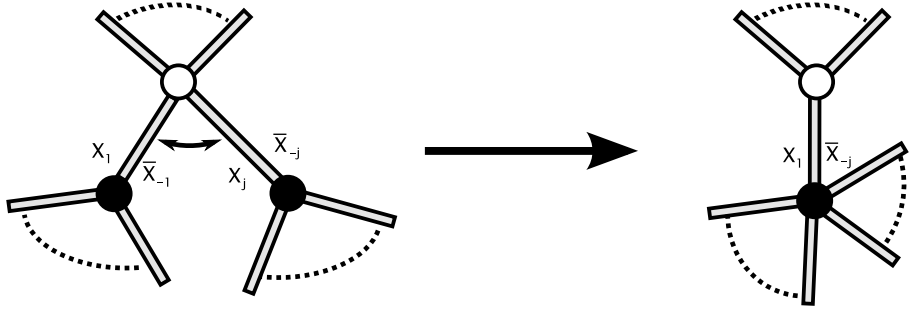}
\caption{ Here two black vertices are connected together through
untwisted edges adjacent at a white vertex.  This bipartite M\"obius
graph reduces to one with one less vertex and one less edge.  The
reduction is found by gluing the two black vertices together and
gluing the two ribbons together along their adjacent sides, here
labeled by $\bar{X}_{-1}$ and $X_j$.  The same reduction would apply
if the two edges were twisted as we could pass this twist through
the white vertex. \label{two-gluing-orient}}
\end{figure}

\begin{figure}
\includegraphics[width=12cm]{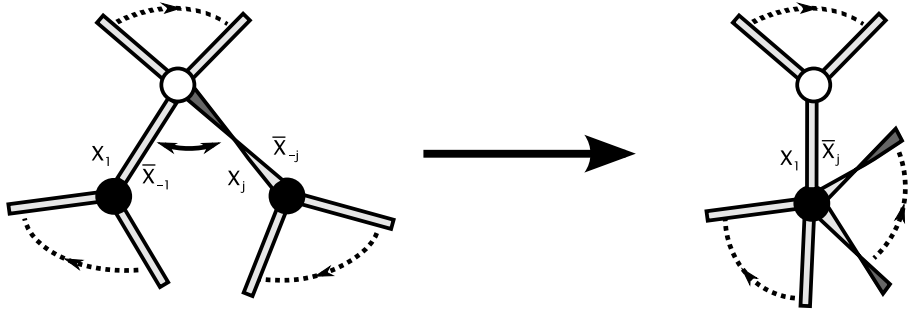}
\caption{ Here two black vertices are connected together through one
twisted edge and one untwisted edge adjacent at a white vertex.
This bipartite M\"obius graph reduces to one with one less vertex
and one less edge.  The reduction is found by gluing the two black
vertices together, gluing the two ribbons together along the sides
adjacent at the white vertex, and reversing both the order and
orientations of remaining ribbons on the second black vertex.
\label{two-gluing-unorient}}
\end{figure}

Therefore we will only need to prove the result for the single black vertex case of
the induction step.
We wish to show that
\begin{equation}
\E( \bar{X}_{-1} X_1 \bar{X}_{-2} X_2 \cdots \bar{X}_{-n} X_n ) = 4^{n-1}
(-2)^{\chi(\Gamma)}\;,
\end{equation}
where $\Gamma$ is a bipartite M\"obius graph with a single black vertex and
half ribbons labeled by $X_{\pm k}$.  We will
do this by induction, there are two cases:
\begin{enumerate}

\item[\bf Case 1:] $X_{-1} = X_j$ for $1 \leq j \leq n$,
\begin{align*}
\E( \bar{X}_{-1} X_1 \cdots \bar{X}_{-j} X_j \cdots X_n ) &=
\E( \bar{X}_{-1} X_j ) \E( \Re( X_1 \cdots \bar{X}_{-j} ) \Re( \bar{X}_{-j-1}
X_{j+1} \cdots X_n)) \\
&= 4 \E(\Re( \bar{X}_{-j} X_1 \cdots X_{j-1} ) \Re( \bar{X}_{-j-1} X_{j+1}
\cdots X_n) )  \,.
\end{align*}
This corresponds to the reduction of the bipartite M\"obius graph pictured in Figure
\ref{reduce-1}, which for $j>1$ splits the single black vertex into two black vertices,
and glues the two edges labeled as $\bar{X}_{-1}$ and $X_j$ together.

The edges $(1,-1,j,-j)$ are adjacent at the white vertex and appear either in this order, or in the reverse order. Thus after the removal of $\{-1,j\}$ we get an ordered pair of labels $(1,-j)$ or $(-j,1)$ to glue back into a ribbon.
On the black vertex, due to our conventions the edges of the ribbons appear in the following order
 $$((-1,1), (-2, 2),\dots, (-j,j), (-k_1,k_1),\dots,(-k_r,k_r)).$$
 Once we split the black vertex into two vertices with the edges of ribbons given by
$((-1,1), (-2, 2),\dots, (-j,j))$ and $((-k_1,k_1),\dots,(-k_r,k_r))$, the removal of $\{-1,j\}$  creates a new pair $(1,-j)$ which we use to create the ribbon to the white vertex.
[After this step, we relabel all edges to use again $\pm 1,\pm 2,\dots$  consecutively.]

We note that the number of faces of the new graph is the same as the previous graph -- the face with the sequence of edges $$(\dots, X_{k_r},\bar{X}_{-1},X_j,\bar{X}_{-k_1},\dots)$$ becomes the face with edges $(\dots, X_{k_r},\bar{X}_{-k_1},\dots)$ on the new graph.
 The
Euler characteristic becomes $\chi_2 = (v_1 + 1) - (e_1 - 1) + f_1 = \chi_1 +
2$ where $v_1, e_1$ and $f_1$ are the number of vertices, edges, and faces of
$\Gamma$.
By the induction assumption we then find
\begin{equation*}
\E( \bar{X}_{-1} X_1 \cdots X_n) = 4 \left[ 4^{(n-1) - 2} (-2)^{\chi_2}
\right] = 4^{n-2} (-2)^{\chi_1 + 2} = 4^{n-1} (-2)^{\chi_1} \,.
\end{equation*}

\begin{figure}
\includegraphics[width=9cm]{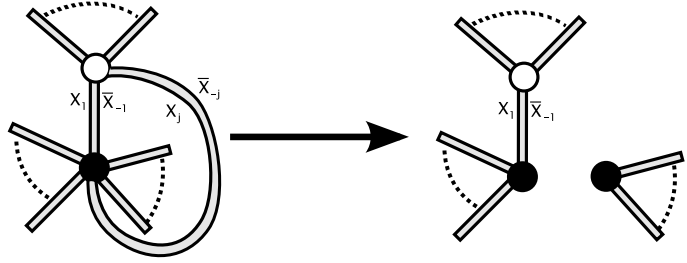}
\caption{ Here we have a black vertex with two ribbons, both twisted
or both untwisted, adjacent at the same white vertex.  The reduction
glues these two ribbons together along their common side.  The
result is a bipartite M\"obius graph with one more vertex and one
less edge.  The resulting graph may or may not be disconnected at
this point. \label{reduce-1}}
\end{figure}

\item[\bf Case 2:] $X_{-1} = X_{-j}$ for $1 < j \leq n$,
\begin{multline*}
\E( \bar{X}_{-1} X_1 \cdots \bar{X}_{-j} X_j \cdots X_n ) \\=
\E( \bar{X}_{-1} \bar{X}_{-j} ) \E( \bar{X}_{j-1} X_{-j+1} \cdots \bar{X}_1
X_j \bar{X}_{-j-1} X_{j+1} \cdots X_n) \\
= (-2)  \E( \bar{X}_{j-1} X_{-j+1} \cdots \bar{X}_1
X_j \bar{X}_{-j-1} X_{j+1} \cdots X_n)  \,.
\end{multline*}
This corresponds to the reduction of the M\"obius graph pictured in Figure
\ref{reduce-2}, which switches the order and the orientations of the
ribbons on one side of the
black vertex from the $\pm 1$ and $\pm j$ ribbons and glues the two
edges adjacent at the white vertex together as shown.
As previously, the removed edges are adjacent at the white vertex. At the black vertex, the labeled lines for the construction of the bipartite graph change from the sequence
$$
(-1,1),(-2,2),\dots,(-j+1,j-1),(-j,j),(-k_1,k_1),\dots,(-k_r,k_r)
$$
to the sequence
$$
(j-1,-j+1),\dots,(2,-2),(1,j),(-k_1,k_1),\dots,(-k_r,k_r)
$$
which then needs to be relabeled to use $\pm 1,\dots,\pm n$. Again the number of faces on the bipartite graphs is preserved:
the face with edges
$$(\dots, 2, -1, -j,k_1,\dots)$$
becomes the face
$$
(\dots,2,1,j,k_1,\dots).
$$
The Euler
characteristic becomes $\chi_2 = v_1 - (e_1 -1) + f_1 = \chi_1 + 1$.  By the
induction assumption we then find
\begin{equation*}
\E( \bar{X}_{-1} X_1 \cdots X_n) = (-2) \left[ 4^{(n-1) - 1} (-2)^{\chi_2}
\right] = (-2) 4^{n-2} (-2)^{\chi_1 + 1} = 4^{n-1} (-2)^{\chi_1} \,.
\end{equation*}

\begin{figure}
\includegraphics[width=9cm]{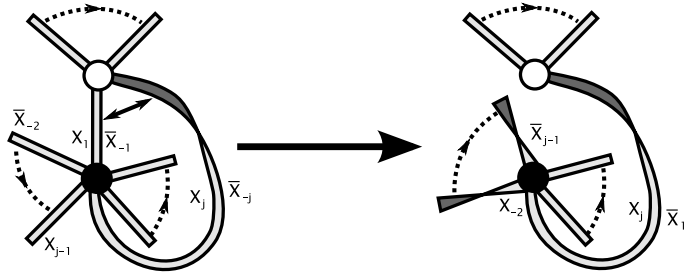}
\caption{ Here we have a black vertex with two ribbons, one twisted
and the other untwisted, adjacent at the same white vertex.  The
reduction glues these two ribbons together along their common side.
The result is a bipartite M\"obius graph with one less edge, and
with the ribbons on one side of the removed ribbon now with reversed
order and orientations. \label{reduce-2}}
\end{figure}

\end{enumerate}

With these two cases checked, by the induction hypothesis, the proof is
completed.

\end{proof}

One should note that this is fundamentally the same proof as in the Wigner
case, however in this case the geometric reduction is given by gluing together
two ribbons, while in the Wigner case the geometric reduction is the
elimination of one ribbon at a time.  The inductive steps remain
the same.

\section{\label{duality} Duality between real and symplectic ensembles}
By $\calM_{M\times N}(\HH )$ we denote the set of all $M\times N$ matrices with entries from $\HH$.
For $\mA\in \calM_{M\times N}(\HH)$, the adjoint matrix is $A^*_{i,j}:=\overline{A}_{j,i}$. The trace is
$\tr(\mA)=
    \sum_{j=1}^NA_{jj}
$.
Since the traces $\tr(\mA)$ may fail to commute, in the formulas we will use $\Re(\tr(\mA))$, compare to
 \cite{Hanlon-Stanley-Stembridge-92}.

\subsection{Duality between GOE and GSE ensembles}

The Gaussian orthogonal ensembles consist of square symmetric matrices, $\mZ$,
whose independent entries are independent (real) Gaussian random
variables; the off diagonal entries have variance $1$ while the
diagonal entries have variance
$2$.
One may show  that
\begin{theoremA} \label{thm3.1} \tolerance=2000
For $\mZ$ from the $N\times N$ Gaussian orthogonal ensemble:
\begin{equation*}
\frac{1}{N^{n-m}} \E( \tr(\mZ^{j_1}) \tr(\mZ^{j_2}) \dots
\tr(\mZ^{j_m}) ) = \sum_{\Gamma} N^{\chi(\Gamma)}  \,,
\end{equation*}
where the sum is over labeled M\"obius graphs $\Gamma$ with $m$
vertices of degree $j_1, j_2, \dots, j_m$, $\chi(\Gamma)$ is the Euler
characteristic and $j_1 + j_2 + \dots + j_m = 2n$.
More generally, if $\mZ_1, \dots, \mZ_s$ are independent $N\times N$ GOE
ensembles and $t: \{ 1, 2, \dots, 2n\} \to \{ 1, \dots, s\}$ is fixed,
then
\begin{multline*}
\frac{1}{N^{n-m}} \E( \tr( \mZ_{t(1)} \dots \mZ_{t(\beta_1)} )
\tr( \mZ_{t(\alpha_2)} \dots \mZ_{t(\beta_2)} ) \times \dots \\ \dots\times
\tr( \mZ_{t(\alpha_m)} \dots \mZ_{t(\beta_m)}) )
= \sum_{\Gamma} N^{\chi(\Gamma)} \,,
\end{multline*}
where $\alpha_1 = 1$, $\alpha_k = j_1 + j_2 +  \dots + 
 j_{k-1} + 1 $, and
$\beta_k = j_1 + j_2 + \dots + j_k $ denote the ranges under the
traces, and where the sum is over labeled color-preserving M\"obius
graphs $\Gamma$ with vertices of degree $j_1, j_2,\dots, j_m$ whose edges are
colored by the mapping $t$.
If there are no $\Gamma$ that are consistent with the coloring we
interpret the sum as being $0$.
\end{theoremA}
The single color version of this Theorem was given in Ref.
\cite{Goulden-Jackson-97}.
\tolerance=1000

The  Gaussian symplectic ensembles (GSE)
consist of square self-adjoint matrices
\begin{equation}
  \label{GUE}  \mZ=\left[Z_{i,j}\right],
\end{equation}
where $\{Z_{i,j}:i<j\}$ is a family of independent $\HH$-Gaussian random variables \eqref{HH}, $\{Z_{i,i}\}$ is a family
of independent real Gaussian random variables of variance $2$,
and $Z_{i,j}=\bar{Z}_{j,i}$. The law of such a matrix has a density
$C\exp(-\frac14 \Re\tr(\mx^2))$,
which is supported on the self-adjoint subset of
$\calM_{N\times N}(\HH)$. ($C=C(N)$ is a normalizing constant.)



In this section we will demonstrate a multi-matrix version of the theorem of Mulase and
Waldron \cite{Mulase-Waldron-03} from our Wick formula.
This represents an improvement to the argument
of Mulase and Waldron as we will not need to rely on labelings of
the vertices by the quaternions $1, i, j,$ and $k$.  This part of
their argument is encoded in relations (\ref{q1}-\ref{q4}).

We will compute here the expected values of the real traces of
powers of a quaternionic self-dual matrix in the  Gaussian
symplectic ensemble \eqref{GUE} where the off-diagonal matrix entries are
\eqref{HH}.

The basic theorem we will prove is
\begin{theorem}\label{T-MW}
For $\mZ$ from the $N\times N$ Gaussian symplectic ensemble:
\begin{equation}
  \label{MW1}
  \frac{1}{(4N)^{n-m}} \E( \Re( \tr( \mZ^{j_1}))
\Re(\tr(\mZ^{j_2})) \dots \Re(\tr(\mZ^{j_m})) ) =
\sum_{\Gamma} (-2N )^{\chi(\Gamma)}\,,
\end{equation}\tolerance=2000
where the sum is over labeled M\"obius graphs $\Gamma$ with $m$ vertices
of degree
$j_1, j_2, \dots, j_m$,  $\chi(\Gamma)$ is the Euler
characteristic and $j_1 + j_2 + \dots + j_m = 2 n$.

More generally, suppose $\mZ_1,\dots,\mZ_s$ are independent $N\times N$  $GSE$ ensembles and $t:\{1,\dots,2n\}\to \{1,\dots, s\}$ is fixed.
Let $\alpha_1=1$, $\alpha_k=j_1 + j_2 + \dots + j_{k-1}+1$, $\beta_k=j_1 + j_2 + \dots + j_{k}$ denote the ranges under traces.
Then
\begin{multline}
  \label{MW expansion}
  \frac{1}{(4N)^{n-m}} \E\big( \Re( \tr( \mZ_{t(1)}\dots \mZ_{t(\beta_1)}))
\Re(\tr(\mZ_{t(\alpha_2)}\dots \mZ_{t(\beta_2)}))\times \dots \\\dots\times \Re(\tr(\mZ_{t(\alpha_m)}\dots \mZ_{t(\beta_m)} ))\big) =
\sum_{\Gamma} (-2N )^{\chi(\Gamma)}\,,
\end{multline}
where the sum is over labeled color-preserving M\"obius graphs $\Gamma$ with vertices
of degree
$j_1$, $j_2$, $\dots$, $j_m$ that are colored with $s$ colors by the mapping $t$. As previously,
  $\chi(\Gamma)$ is the Euler
characteristic and $j_1 + j_2 + \dots + j_m = 2 n$. (If there are no $\Gamma$ that are consistent with the
 coloring $t$, we interpret the sum as $0$.)
\end{theorem}
\tolerance=1000

In view of Theorem \ref{thm3.1}, Theorem \ref{T-MW} gives the
duality between the GOE and GSE ensembles of random matrices as $N
\to -2N$.

\begin{example}
  To illustrate the multi-matrix aspect of the theorem, suppose $\mZ_1,\mZ_2$ are independent
  $N\times N$ GSE  ensembles and $\mX_1,\mX_2$
are  independent
  $N\times N$  GOE ensembles. For $p,r,k\geq 0$
   let $b_{p,r,k}(N):=\E((\tr(\mX_1^{2p}\mX_2^{2r}))^k)$.
  By Theorem \ref{thm3.1}, $b_{p,r,k}(N)=N^{(p +r-1)k}\sum_{\Gamma} N^{\chi(\Gamma)}$.
  Therefore Theorem \ref{T-MW} implies that the moments for the independent
  GSE ensembles are determined from the corresponding moments of independent GOE ensembles
  by the ``dual formula''
$$\E((\tr(\mZ_1^{2p}\mZ_2^{2r}))^k)=(-2)^{(p +r-1)k}b_{p,r,k}(-2N).$$
In particular,
$$\E(\tr(\mZ_1^{2p}))=(-2)^{p-1} b_{p,0,1}(N).$$
This is equivalent to \cite[Theorem 6]{Ledoux-07} with $c_p(N)=\E(\tr(\mZ^{2p}))/2^p$;
an extra factor of $2^p$ in our formula is due to the fact that
the symplectic ensemble in \cite{Ledoux-07}
has the variance of $1/2$ instead of our choice of $1$.
\end{example}

\begin{proof}[Proof of Theorem \ref{T-MW}]

We begin by expanding out the traces in terms of the matrix entries
\begin{align}\label{star1}
 &\frac{N^m}{4^{n-m}} \E(\Re(\tr(\mZ^{j_1}))
 \Re(\tr(\mZ^{j_2})) \dots \Re(\tr(\mZ^{j_m})) ) =\\
 =& \nonumber
 \sum_{\begin{matrix}
 1\leq a_1, a_2, \dots, a_{j_1} \leq N \\
 1\leq b_1, b_2, \dots, b_{j_2} \leq N \\ \vdots \\
 1 \leq c_1, c_2, \dots, c_{j_m} \leq N
 \end{matrix} }
 \begin{matrix}
  \frac{1}{(4n)^{n-m}} \E(\Re( Z_{a_1, a_2} Z_{a_2, a_3}
 \dots Z_{a_{j_1}, a_1}) \Re( Z_{ b_{1}, b_{2} } \dots
 Z_{b_{j_2}, b_1} ) \times\dots
\\
\hfill \dots \times \Re( Z_{ c_{1}, c_{2} } \dots
 Z_{ c_{j_m}, c_1} ) )\,.
 \end{matrix}
 \end{align}
Note that essentially the same expansion applies to \eqref{MW expansion}, except that the consecutive
entries $Z_{i,j}^{(t)}$ must now be labeled
also by the ``color" $t$.
 It is then better to index the products
by the cycles of a permutation. Put 
$$\sigma=(1,\dots,\beta_1)(\alpha_2,\dots,\beta_2)\ldots(\alpha_m,\dots ,2n)\,.$$
In this notation, the expansion of the  left hand side of \eqref{MW expansion} takes the form
 $$\sum_{a:\{1\dots2n\}\to\{1\dots N\}}\E\left(\prod_{c\in\sigma}
 \Re(\prod_{j\in c}[\mZ_{t(j)}]_{a(j),a(\sigma(j))})\right).$$
(Since  $\Re$   has the cyclic property
\ref{tmp**}, this expression is well defined.)

From \eqref{Wick1} it follows that the right hand side of
\eqref{star1} can be expanded as a sum over all pairings, and we can
assume that the pairs are independent.

Of course, pairings that match two independent random variables do not contribute to the sum.
The pairings that contribute to the sum are of three different types: pairs that match $Z$ with another $Z$ at a different position in the product,
pairs that match $Z$ with $\bar{Z}$, and pairings that match the diagonal (real) entries.

We first dispense with the pairings that match a diagonal entry $Z_{i,i}=\xi$ with another diagonal entry
$Z_{j,j}=\xi$. Since real numbers commute with quaternions,
$$\Re(\E(q_0\xi q_1\xi q_2))=\E(\xi^2) \Re(q_0q_1q_2)=2 \Re(q_0q_1q_2).$$
On the other hand,  using the cyclic property of $\Re(\cdot)$ and adding formulas \eqref{q1} and \eqref{q2}, we see that the same answer arises from
$$\E\Re(q_0Z q_1Z q_2)+\E\Re(q_0Z q_1\bar{Z} q_2)=\E\Re(Z q_1Z q_2q_0)+\E\Re(Z q_1\bar{Z} q_2q_0)=2\Re(q_1q_2q_0).$$
So the contribution of each such diagonal pairing is the same as that of two matches of quaternionic entries $(Z,Z)$ and $(Z,\bar{Z})$.

Thus, once we replace all the real entries that came from the diagonal entries by the corresponding quaternion-Gaussian pairs,
we get the sum over all possible pairs of matches of the first two types only.
We label all such pairings by M\"obius graphs, with the interpretation that
pairings or random variables $(Z,Z)$ correspond to twisted ribbons, while pairings $(Z,\bar{Z})$
correspond to ribbons without a twist. The pairs of variables at different ribbons can now be assumed independent.
In the multi-matrix case, the edges of the graph at each vertex are colored according to the function $t$,
which restricts the number of available pairings.

We now relabel the $Z_{i,j}$ by $Z_{a_k, a_{k+1}} = X_k$, $Z_{b_k, b_{k+1}}
= X_{j_1 + k} $, $\dots$,
$Z_{c_k, c_{k+1}} = X_{j_1 + j_2 + \dots
j_{m-1} + k} $. Our claim is then that given a M\"obius graph $\Gamma$
with vertices of degree $j_1$, $j_2, \dots,$ and $j_m$ with edges
labeled by $X_k$, satisfies
\begin{multline}
\frac{1}{(4N)^{n-m}} \E( \Re(X_1 \dots X_{j_1}) \Re( X_{j_1+1}
\dots X_{j_1+j_2} ) \dots \Re( X_{j_1 + \dots +j_{m-1}+1} \dots
X_{2n}) )\\= (-2N)^{\chi(\Gamma)} N^{-f(\Gamma)} \,,
\end{multline}
where $f(\Gamma)$ is the number of faces of $\Gamma$.
The $N^{-f(\Gamma)}$ terms are removed by the summations in
(\ref{star1}), the relations given by the edges of $\Gamma$ reduce
the number of summations leaving us with $f(\Gamma)$ sums from $1$
to $N$.  Collecting the powers of $N$ we find
that this is the same as \eqref{star2}.
\end{proof}

\subsection{Duality between the Laguerre Ensembles}

\begin{figure}
\includegraphics[width=8cm]{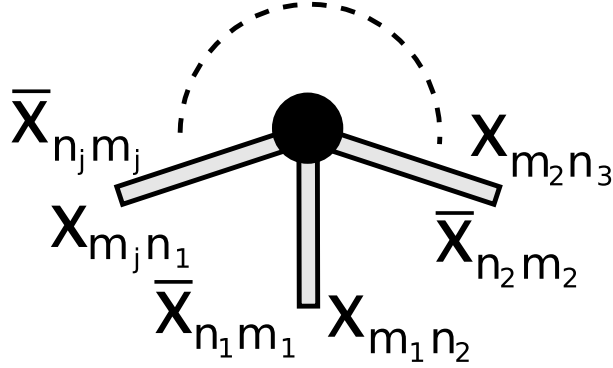}
\caption{ For the Wishart applications we will label the sides of
the ribbons emanating from a black vertex.  One side receives the
label $X_{n_k m_k}$ while the other receives the label $\bar{X}_{m_k
n_{k+1} }$.  Note that the interior index, $m_k$, of the two labels
on each ribbon match; while the exterior index, $n_k$, of the labels
on ribbons adjacent at a black vertex match. \label{black-white-3}}
\end{figure}

We will now consider the analogous problem for quaternionic Wishart random
matrices. These are $N\times N$ matrices of the form
\begin{equation}
  \label{H-Wishart}
  \mW = \mX^* \mX\,,
\end{equation}
where $\mX$ is an $M\times N$ matrix with independent $\HH$-Gaussian entries \eqref{HH}.
The Wick formula (\ref{Wick1}) in this case is expressed in terms of
bipartite M\"obius graphs defined in section \ref{Sect bipartite};
although in this case it is helpful to view our labeling of the ribbons around
each black vertex in a slightly different way (see Figure
\ref{black-white-3}).
Our heuristic is the following: when edges are glued to a white
vertex we identify the labels on the interior of each ribbon as all
equal, while the labels on the outside of each ribbon are equal to
the others which border the same face.  The reason for this
dichotomy is that $\mX$ is not a square matrix and therefore the
$n_j$'s and $m_j$'s are to be treated as separate families of
variables.

In the multi-matrix case, we color the ribbons by the number of a matrix, and color the white vertices so that only the ribbons of matching color can be attached.

\begin{figure}
\begin{center}\includegraphics[width=5cm]{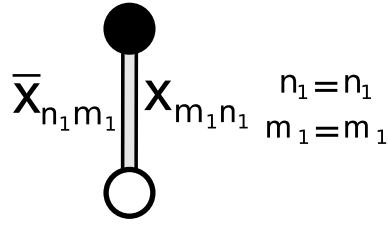}
\end{center}
\caption{ The only possible bipartite M\"obius graph with a single
black vertex of degree one.  In this case the graph is embedded on a
copy of the Riemann sphere. \label{black-white-ex1}}
\end{figure}

\begin{figure}
\begin{center}\includegraphics[width=10cm]{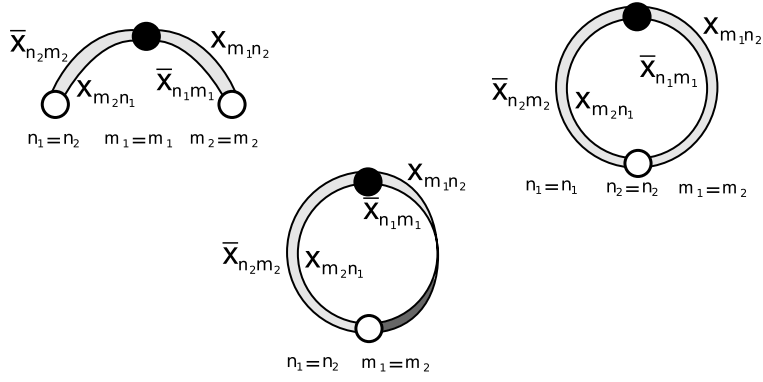}\end{center}
\caption{ The three possible bipartite M\"obius graphs with a single
black vertex of degree two.  In the top two cases the graph is
embedded on a copy of the Riemann sphere, while in the bottom case
the graph is embedded on a copy of the projective sphere.
\label{black-white-ex2}}
\end{figure}

For example,
\begin{equation*}
\E( \tr(\mW) ) = \sum_{m=1}^M \sum_{n=1}^N \E( X^*_{nm} X_{mn} )
= 4 M N
\end{equation*}
corresponds to the bipartite M\"obius graph in Figure \ref{black-white-ex1}.
As another example
\begin{align*}
\E(\tr(\mW^2)) &= \sum_{m_1, m_2 =1}^M \sum_{n_1, n_2 =1}^N
\E( X^*_{n_1 m_1} X_{m_1 n_2} X^*_{n_2 m_2} X_{m_2 n_1} )
\\&= \sum_{m_1, m_2=1}^M \sum_{n_1, n_2=1}^N\Big[
\E( X^*_{n_1 m_1} X_{m_1 n_2} ) \E( X^*_{n_2 m_2} X_{m_2 n_1})
\\&\phantom{\sum_{m_1, m_2=1}^M \sum_{n_1, n_2=1}^N }
+ \E( X^*_{n_1 m_1} X^*_{n_2 m_2} ) \E( \bar{X}_{m_1 n_2} X_{m_2 n_1})
\\&\phantom{\sum_{m_1, m_2=1}^M \sum_{n_1, n_2=1}^N }
+ \E( X^*_{n_1 m_1} X_{m_2 n_1}) \E( X_{m_1 n_2} X^*_{n_2 m_2} )\Big]
\\&= 16 N M^2 - 8 N M + 16 N^2 M\,,
\end{align*}
corresponds to the bipartite M\"obius graphs in Figure \ref{black-white-ex2}.

In this case our theorem takes the form.
\begin{theorem}
With $M = \lambda N$ we find
\[ \E( \tr( \mW^{j_1} ) \tr( \mW^{j_2}) \cdots \tr( \mW^{j_m}) )
= N^{n-m}   4^{n-m} \sum_{\Gamma} \lambda^{w(\Gamma)} (-2N )^{\chi(\Gamma)}\;,
\]
where the sum is over the bipartite M\"obius graphs with black
vertices having degrees $j_1, j_2, \dots, j_m$, $n = j_1 + j_2 +
\dots + j_m $ is the number of edges, $w(\Gamma)$ is the number of
white vertices, and $\chi(\Gamma)$ is the Euler characteristic.

More generally, suppose that $\mW_1,\mW_2,\dots,\mW_s$ are $N\times N$ independent quaternionic Wishart with parameters $M_1,\dots,M_s$. Denote $\la_j=M_j/N$.
Fix  $t:\{1,\dots,n\}\to \{1,\dots, s\}$, and
let $\alpha_1=1$, $\alpha_k=j_1 + j_2 + \dots + j_{k-1}+1$, $\beta_k=j_1 + j_2 + \dots + j_{k}$ denote the ranges under traces.
Then
\begin{multline}
  \label{Wishart expansion}
  \frac{1}{(4N)^{n-m}} \E\big( \Re( \tr( \mW_{t(1)}\dots \mW_{t(\beta_1)}))
\Re(\tr(\mW_{t(\alpha_2)}\dots \mW_{t(\beta_2)})) \times\dots
\\\dots\times  \Re(\tr(\mW_{t(\alpha_m)}\dots \mW_{t(\beta_m)} ))\big)=
\sum_{\Gamma} \lambda_1^{w_1(\Gamma)}\lambda_2^{w_2(\Gamma)}\dots \lambda_s^{w_s(\Gamma)} (-2N )^{\chi(\Gamma)},
\end{multline}
where the sum is over the bipartite M\"obius graphs with $m$ black vertices of degrees $j_1,\dots,j_m$, whose ribbons are colored according to $t$ by one of the
colors $1,\dots,s$,
 and with white vertices that  are colored to match the ribbons. Here $w_j(\Gamma)$ is the number of white vertices of color $j$,   $\chi(\Gamma)$ is the Euler
characteristic and $j_1 + j_2 + \dots + j_m = n$ is the number of edges/ribbons.
\end{theorem}

\begin{proof}
We begin by expanding out the traces in terms of the matrix entries of
$\mX$ where $\mW = \mX^* \mX $ and $\mX$ is an $M \times N$ matrix of
quaternionic Gaussian random variables
\begin{align} \label{Wishart-expansion}
& \frac{N^{m-n}}{4^{n-m}} \E( \tr( \mW^{j_1}) \tr( \mW^{j_2})
\cdots \tr( \mW^{j_m}) )  = \\ \nonumber
& \sum_{
\begin{matrix}
 1 \leq a_1, a_2, \dots, a_{j_1} \leq N \\
1 \leq A_1, A_2, \dots, A_{j_1} \leq M \\
1 \leq b_1, b_2, \dots, b_{j_2} \leq N \\
1 \leq B_1, B_2, \dots, B_{j_2} \leq M \\
\vdots\\
1 \leq C_1, C_2, \dots, C_{j_m} \leq M \end{matrix}}
\begin{matrix}
\frac{N^{m-n}}{4^{n-m}} \E( \Re( [\mX^*]_{a_1 A_1} [\mX]_{A_1 a_2} [\mX^*]_{a_2
  A_2} [\mX]_{A_2 a_3} \dots [\mX^*]_{a_{j_1} A_{j_1}} [\mX]_{A_{j_1} a_1} )
\\
\hfill \Re( [\mX^*]_{b_1 B_1}  \dots [\mX]_{B_{j_2} b_1} ) \dots
\Re( [\mX^*]_{c_1 C_1} \dots [\mX]_{C_{j_m} c_1} ) )\,.
\end{matrix}
\end{align}
Note that a similar expansion would be found for a multi matrix
expression, see below.

From (\ref{Wick1}) it follows that
(\ref{Wishart-expansion}) can be expanded as a sum over all pairings,
and we can assume that the pairs are independent.

We label all such pairings by bipartite M\"obius graphs, that is
M\"obius graphs with labeled black vertices of degree $j_1, j_2, \dots,
j_m$, and unlabeled white vertices of any degree.  The half edges of
each graph are labeled by $\bar{X}_{A_1 a_1}= \bar{X}_{-1}$, $X_{A_1
  a_2} = X_1 $ and so on.
The relations given by this M\"obius graph reduce the number of sums
over capital indices (which go from 1 to $M$) to $w(\Gamma)$ the
number of white vertices, and reduce the number of sums over lower
case indices (which go from 1 to $N$) to $f(\Gamma)$ the number of
faces of $\Gamma$.  Therefore from Theorem \ref{thm2.1},
with $M = \lambda N$, the total
contribution from the Wishart M\"obius graph $\Gamma$ to the sum is
\[ \frac{N^{m-n}}{4^{n-m}}
\left[ 4^{n-m}  \lambda^{w(\Gamma)} (-2)^{\chi(\Gamma)} \right]
N^{w(\Gamma) + f(\Gamma)}
= \lambda^{w(\Gamma)}  (-2N)^{\chi(\Gamma)}\,, \]
as $\chi(\Gamma) = m + w(\Gamma) - n + f(\Gamma)$.

The multimatrix version of this proof requires no conceptual changes,  but  the notation becomes more cumbersome.
 In this case it is more convenient to index the products
by the cycles of a permutation. Put 
$$\sigma=(1,\dots,\beta_1)(\alpha_2,\dots,\beta_2)\ldots(\alpha_m,\dots, n).$$ 
In this notation, the expansion on the right hand side of \eqref{Wishart-expansion} is replaced by
\begin{multline*}
\sum_{a:\{1\dots n\}\to\{1\dots N\}}\E\left(\prod_{c\in\sigma}\Re(\prod_{j\in c}
[\mW_{t(j)}]_{a(j),a(\sigma(j))}
)
\right)\\
=\sum_{a:\{1\dots n\}\to\{1\dots N\}}\sum_{b\in\calB(t)} \E\left(\prod_{c\in\sigma}\Re(\prod_{j\in c}
\overline{[{\mX}_{t(j)}]}_{b(j)a(j)}[\mX_{t(j)}]_{b(j),a(\sigma(j))})\right),
\end{multline*}
where $\calB(t)
=\{1,\dots,M_{t(1)}\}\times\{1,\dots,M_{t(2)}\}\times\dots\times \{1,\dots,M_{t(n)}\}$.

From (\ref{Wick1}) the expected value can again be expanded as a sum over the pairings,
where we can assume that the pairs are independent random variables and the pairings connect only the random variables of the same color $t(j)$. For each bipartite M\"obius graph $\Gamma$,
 the sum over $b\in\calB(t)$ contributes a factor of $M_j^{w_j(\Gamma)}$ per color, where $w_j(\Gamma)$ is the number of vertices of color $j$ in $\Gamma$, while the sum over $a$ contributes $N^{f(\Gamma)}$, as in the one-matrix case.
 Therefore from Theorem \ref{thm2.1}, the left hand side of \eqref{Wishart expansion} becomes
\begin{multline*}
(4N)^{m-n}\sum_\Gamma   M_1^{w_1(\Gamma)}\dots M_s^{w_s(\Gamma)}N^{f(\Gamma)} 4^{n-m}(-2)^{\chi(\Gamma)}\\
 =
 N^{m-n}\sum_\Gamma  N^{w(\Gamma)+f(\Gamma)}(-2)^{\chi(\Gamma)}\prod_{j=1}^s \la_j^{w_j(\Gamma)}\\
 =
\sum_\Gamma  N^{w(\Gamma)+m+f(\Gamma)-n} (-2)^{\chi(\Gamma)}\prod_{j=1}^s \la_j^{w_j(\Gamma)}\\=
\sum_\Gamma  N^{v(\Gamma)+f(\Gamma) 
-e(\Gamma)} (-2)^{\chi(\Gamma)}\prod_{j=1}^s \la_j^{w_j(\Gamma)}=\sum_\Gamma  (-2N)^{\chi(\Gamma)}\prod_{j=1}^s \la_j^{w_j(\Gamma)} .
\end{multline*}
\end{proof}
We remark that the result is again in duality to the
formula for real Wishart matrices, where the analogue of \eqref{Wishart expansion} has the
following form:
\begin{multline}
  \label{Real Wishart expansion}
  \frac{1}{N^{n-m}} \E\big(  \tr( \mW_{t(1)}\dots \mW_{t(\beta_1)})
\tr(\mW_{t(\alpha_2)}\dots \mW_{t(\beta_2)})\times\dots \\ \dots \times \tr(\mW_{t(\alpha_m)}\dots \mW_{t(\beta_m)} )\big)
=
\sum_{\Gamma} \lambda_1^{w_1(\Gamma)}\lambda_2^{w_2(\Gamma)}\dots \lambda_s^{w_s(\Gamma)} N ^{\chi(\Gamma)}.
\end{multline}
For the one-matrix case this can be read out from \cite{Hanlon-Stanley-Stembridge-92}.
For the multivariate case this is a reinterpretation of \cite[Theorem 2.10]{Bryc-07}. 

The following illustrates a one-matrix version of this duality.
\begin{example}
Suppose $\mW_\RR$ is $N\times N$  real Wishart matrix with parameter $M$, and let
$$Q_{j_1,\dots,j_m}^{\RR}(M,N)=\E( \tr( \mW_\RR^{j_1}) \tr( \mW_\RR^{j_2})
\cdots \tr( \mW_\RR^{j_m}) ).$$
If  $\mW_\HH$ is $N\times N$  quaternion Wishart matrix with parameter $M$, and $n=j_1+\dots +j_m$,
then \eqref{Real Wishart expansion} and \eqref{Wishart expansion} give
$$\E( \tr( \mW_\HH^{j_1}) \tr( \mW_\HH^{j_2})
\cdots \tr( \mW_\HH^{j_m}) )=(-2)^{n-m}Q_{j_1,\dots,j_m}^{\RR}(-2M,-2N).$$
This is equivalent to \cite[Corollary 4.2]{Hanlon-Stanley-Stembridge-92}, as with $\la=(j_1,j_2,\dots,j_m)$ in their notation
\begin{multline*}
 Q_{j_1,\dots,j_m}^{\HH}(M,N)=2^{-n}\E( \tr( \mW_\HH^{j_1}) \tr( \mW_\HH^{j_2})
\cdots \tr( \mW_\HH^{j_m}) )\\=(-1)^{n-m} 2^{-m}Q_{j_1,\dots,j_m}^{\RR}(-2M,-2N);
\end{multline*}
an extra factor  of $2^n$ in our formula is due to the fact that the quaternion Gaussian law in \cite{Hanlon-Stanley-Stembridge-92}
has the variance of $1/2$ instead of our choice of $1$.
\end{example}

\subsection*{Acknowledgements}
The authors thank G\'erard Letac for citation \cite{Isserlis-1918}. We would like to thank the referee for several substantial corrections.

\bibliographystyle{acm}
\bibliography{quaternionic-07}
\end{document}